\documentclass[12pt,a4paper]{article}
\usepackage{amsthm,amssymb,amsmath}
\usepackage[a4paper]{geometry}
\usepackage{hyperref,graphicx,enumerate}
\usepackage{tikz}
\hypersetup{
 pdftitle={Perfect matchings in Claw-free Cubic Graphs},
 pdfauthor={Sang-il Oum}
}
\newtheorem{THM}{Theorem}

\newtheorem{PROP}[THM]{Proposition}
\theoremstyle{definition}
\newcommand\C{\mathcal C}

\title{Perfect Matchings in Claw-free Cubic Graphs}
\author{Sang-il Oum\thanks{{\tt sangil@kaist.edu}} \footnote{
    Supported by  the SRC Program of Korea Science and Engineering Foundation
    (KOSEF)
    grant funded by the Korea government (MEST)
    (No. R11-2007-035-01002-0).
  }\\
  Department of Mathematical Sciences, \\
  KAIST, Daejeon, 305-701,
  Republic of Korea.
}
\date{November 9, 2009}
\begin{document}
\maketitle
\begin{abstract}
  Lov\'asz and Plummer conjectured that
  there exists a fixed positive constant $c$ such that
  every cubic $n$-vertex graph with no cutedge
  has at least $2^{cn}$ perfect matchings.
  Their conjecture has been verified for bipartite graphs by Voorhoeve
  and
  planar graphs by Chudnovsky and Seymour.
  We  prove that every claw-free cubic $n$-vertex graph  with no cutedge
  has more than $2^{n/12}$ perfect matchings, thus verifying the
  conjecture for claw-free graphs.
\end{abstract}
\section{Introduction}
A graph is \emph{claw-free} if it has no induced subgraph
isomorphic to $K_{1,3}$.
A graph is \emph{cubic} if every vertex has exactly three incident
edges.
A well-known classical theorem of Petersen \cite{Petersen1891} states
that every cubic graph with no cutedge has a perfect matching.
Sumner \cite{Sumner1974a} and Las Vergnas \cite{LasVergnas1975}
independently
showed that
every connected claw-free graph with even number of vertices
has a perfect matching.
Both theorems imply that every claw-free cubic graph with no cutedge has at least one
perfect matching.

In 1970s, Lov\'asz and Plummer
conjectured that every cubic graph with no cutedge has exponentially
many perfect matchings; see \cite[Conjecture 8.1.8]{LP1986}.
The best lower bound has been obtained by 
Esperet, Kardo\v{s}, and Kr\'al' \cite{EKK2009}.
They showed that the number of perfect matchings
in a sufficiently large cubic graph with no cutedge
always exceeds any fixed linear function in the number of vertices.

So far the conjecture is known to be true for bipartite graphs and
planar graphs.
For bipartite graphs,
Voorhoeve \cite{Voorhoeve1979}  
proved that every \emph{bipartite} cubic
$n$-vertex graph
has at least $6(4/3)^{n/2-3}$ perfect matchings.
Recently, Chudnovsky and Seymour \cite{CS2008}
proved that every \emph{planar} cubic $n$-vertex graph with no cutedge
has at least $2^{n/655978752}$ perfect matchings.

\begin{figure}
  \centering
  \tikzstyle{every node}=[circle, draw, fill=black!50,
  inner sep=0pt,minimum width=4pt]
  \begin{tikzpicture}[thick,scale=0.5]
    \draw \foreach \x in {0,3,6,13}
    {
      (\x,1) node {} --+(0,-2) node {} --+(1,-1) node{} --+ (2,-1)  
      (\x,1) --+ (1,-1)
   }
    \foreach \x in {3,6,13,16}
    {
      (\x,1) node {} --+ (-1,-1) node {} --+ (0,-2) node {} --cycle
    }
    \foreach \x in {-2,18}
    {
      (\x,1) node {} --+ (-1,-1) node {} --+ (0,-2) node {}
      --+(1,-1)      node {} --+ (-1,-1)
      (\x,1) --+ (1,-1)
   }
    \foreach \x in {-2,16}
    {
      (\x,1) --+ (2,0)
      (\x,-1) --+ (2,0)
     }
     (11,0)--(12,0);
    \draw [dashed] (8.5,0)--(10.5,0);
\end{tikzpicture}
  \caption{Claw-free cubic graphs with only 9 perfect matchings}
  \label{fig:ex}
\end{figure}
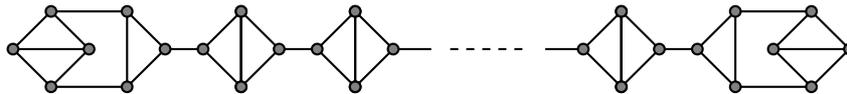

We prove that
every claw-free cubic $n$-vertex graph with no cutedge
has more than
\[2^{n/12}\]
perfect matchings.
The graph should not have any cutedge;
in Figure \ref{fig:ex}, we provide an example of a claw-free cubic graph with only 9 perfect
matchings.

Our approach is to use the structure of 2-edge-connected claw-free cubic graphs.
The \emph{cycle space} $\C(H)$ of~$H$ is a collection
of the edge-disjoint union of cycles of~$H$.
It is well known that $\C(H)$ forms a vector space over $GF(2)$
and \[\dim\C(H)=|E(H)|-|V(H)|+1\] if $H$ is connected, see Diestel \cite{Diestel2005}.
Roughly speaking, almost all 2-edge-connected claw-free cubic graph
$G$
can be built from a
2-edge-connected cubic  multigraph $H$ by certain operations
so that 
every member of $\C(H)$
can be extended to
$2$-factors of $G$.
We will have two cases to consider; either $H$ is big or small. If $H$ is big, then
$\C(H)$ is big enough to prove that $G$ has many
$2$-factors.
If $H$ is small, then
we find a $2$-factor of $H$ using 
many of the specified edges of $H$ so that
when transforming this $2$-factor of $H$ to that of $G$,
each of those edges of $H$ has many  ways to make $2$-factors of $G$.

\section{Structure of 2-edge-connected claw-free cubic graphs}
Graphs in this paper have no parallel edges and no loops,
and multigraphs can have parallel edges and loops.
We assume that a loop is counted twice when measuring a degree of a
vertex in a multigraph.
Every 2-edge-connected cubic multigraph can not have loops
because if it has a loop, then it must have a cutedge.

We describe the structure of claw-free cubic
graphs given by 
Palmer et~al.~\cite{PRR2002}.
A \emph{triangle} of a graph is a set of three pairwise adjacent
vertices.
\emph{Replacing a vertex  $v$ with a triangle} in cubic graph is to
replace $v$ with three vertices $v_1,v_2,v_3$ forming a triangle
so that if $e_1,e_2,e_3$ are three edges incident with $v$,
then $e_1,e_2,e_3$ will be incident with $v_1,v_2,v_3$ respectively.

Every vertex in a claw-free cubic graph is in 1, 2, or 3 triangles.
If a vertex is in 3 triangles, then the component containing the
vertex is isomorphic to $K_4$.
If a vertex is in exactly 2 triangles,
then it  is in an induced subgraph isomorphic to $K_4\setminus e$ for
some edge $e$ of $K_4$.
Such an induced subgraph is called a \emph{diamond}.
It is clear that no two distinct diamonds intersect.

A \emph{string of diamonds}
is a maximal sequence $D_1,D_2,\ldots,D_k$
of diamonds
in which, for each $i\in \{1,2,\ldots,k-1\}$, $D_i$
has a vertex adjacent to a vertex in $D_{i+1}$.
A string of diamonds has exactly two vertices of degree $2$,
which are called the \emph{head} and the \emph{tail} of the string.
\emph{Replacing an edge $e=uv$ with a string of diamonds} with the
head $x$ and the tail $y$
is to remove $e$ and add edges $ux$ and $vy$.

A connected claw-free cubic graph in which every vertex is in a diamond
is called a \emph{ring of diamonds}.
We require that a ring of diamonds contains at least 2 diamonds.
It is now straightforward to describe the structure of
2-edge-connected claw-free cubic graphs as follows.
\begin{PROP}\label{prop}
  A graph $G$
  is 2-edge-connected claw-free cubic 
  if and only if
  either
  \begin{enumerate}[(i)]
  \item $G$ is isomorphic to $K_4$,
  \item $G$ is a ring of diamonds, or
  \item $G$ can be built from a 2-edge-connected cubic  multigraph $H$
    by
    replacing some edges of $H$ with strings of diamonds
    and 
    replacing each vertex of $H$ with a triangle.
\end{enumerate}
\end{PROP}
\begin{proof}
  Let us first prove the ``if'' direction.
  It is easy to see that $G$ is 2-edge-connected cubic
  and has no loops or parallel edges.
  If $G$ is built as in (iii), then
  clearly $G$ has neither loops nor parallel edges,
  and every vertex of $G$ is in a triangle and therefore $G$ is
  claw-free. Note that since $H$ is $2$-edge-connected, $H$ can not
  have loops.
  
  To prove the ``only if'' direction, 
  let us assume that $G$ is a 2-edge-connected claw-free cubic graph.
  We may assume that $G$ is not isomorphic to $K_4$
  or a ring of diamonds.
  We claim that $G$ can be built from a 2-edge-connected cubic  multigraph
  as in (iii). Suppose that $G$ is a counter example with the minimum
  number of vertices.
  
  If $G$ has no diamonds, then
  every vertex of $G$ is in exactly one triangle
  and therefore $V(G)$ can be partitioned into disjoint triangles.
  By contracting each triangle, we obtain a 2-edge-connected cubic
   multigraph $H$.

  So $G$ must have a string of  diamonds.
  Let $D$ be the set of vertices in the string of diamonds.
  Since $G$ is cubic,
  $G$ has two vertices not in $D$, say $u$ and $v$, adjacent to $D$.
  If $u=v$, then because the degree of $u$ is $3$,
  $u$ must have another incident edge $e$ but $e$ will be a cutedge of $G$.
  Thus $u\neq  v$.

  If $u$ and $v$ are adjacent in $G$,
  then $u$ and $v$ must has a common neighbor $x$, because otherwise
  $G$ will have an induced subgraph isomorphic to $K_{1,3}$.
  However one of the edges incident with $x$ will be a cutedge of $G$, a
  contradiction.

  Thus $u$ and $v$ are nonadjacent in $G$.
  Let $G'=(G\setminus D)+uv$, that is obtained from $G$ by
  deleting $D$ and adding an edge $uv$.
  Then $G'$ has no parallel edges or loops
  and moreover $G'$ is 2-edge-connected claw-free cubic.
  Since $G$  has a vertex not in a diamond,
  so does $G'$
  and therefore
  $G'$ can be built from a 2-edge-connected cubic  multigraph $H$
  by replacing some edges with strings of diamonds
  and replacing each vertex of $H$ with a triangle.
  Since $D$ is chosen maximally,
  $u$ and $v$ are not in diamonds
  and therefore $H$ has the edge $uv$.
  So we can obtain $G$ from $H$ by doing all replacements to obtain
  $G'$ and then replacing the edge $uv$ with a string of diamonds.
  This completes the proof.
\end{proof}
We remark that Proposition~\ref{prop} can be seen as a corollary of
the structure theorem of quasi-line graphs by Chudnovsky and Seymour \cite{CS2005}.
A graph is a  \emph{quasi-line} graph  if the neighborhood of each vertex is
expressible as the union of two cliques.
It is obvious that every claw-free cubic graph is a quasi-line graph.
Chudnovsky and Seymour \cite{CS2005} proved that every connected quasi-line graph is
either a fuzzy circular interval graph or a composition of fuzzy
linear interval strips.
For $2$-edge-connected claw-free cubic graphs, a fuzzy circular interval graph
corresponds to a ring of diamonds
and a composition of fuzzy linear interval strips corresponds to the
construction (iii) of Proposition~\ref{prop}.
\section{Main theorem}

\begin{figure}
  \centering
 \begin{tikzpicture}[thick,scale=.8]
  \draw  (11.5cm,2cm) node [] {or};
  \tikzstyle{every node}=[circle, draw, solid,fill=black!50,
  inner sep=0pt,minimum width=4pt]
    \draw [dotted]
    (-90:1) -- (0,0) node{}
    (30:1)-- (0,0)  -- (150:1)
    [xshift=3cm] (30:1)--(30:.5)
    (-90:1)--(-90:.5)
    (150:1)--(150:.5);
    \draw [xshift=3cm] (30:.5) node {}--(150:.5) node {}--(-90:.5)
    node  {}--cycle;
    \draw [dotted,yshift=2cm]    (-90:1) -- (0,0) 
    [xshift=3cm]    (-90:1)--(-90:.5)    (150:.5)--(30:.5)   ;
   \draw  [yshift=2cm]  (30:1) -- (0,0) node {} -- (150:1) 
   [xshift=3cm] (30:1) --(30:.5) node {}--(-90:.5) node
   {}--(150:.5) node {}--(150:1) ;

   \draw [dotted,xshift=6cm]  (0,0) --(2,0)
   [xshift=4.5cm] (0,0)--(.5,0)  (1.5,0)--(2,0)  (1,.5)--(1,-.5);
   \draw [xshift=10.5cm] (.5,0)node{}--(1,.5)node{}--(1.5,0)node{}--(1,-.5)node{}--cycle;
   
   \draw [xshift=6cm,yshift=2cm]  (0,0)--(2,0)
   [xshift=3cm] (0,0)--(.5,0)  (1.5,0)--(2,0)  (1,.5)--(1,-.5)
   [xshift=3cm] (0,0)--(.5,0)  (1.5,0)--(2,0)  (1,.5)--(1,-.5);
   \draw [dotted,xshift=9cm,yshift=2cm] (.5,0)--(1,.5)
   (1.5,0)--(1,-.5)
   [xshift=3cm](.5,0)--(1,-.5)   (1.5,0)--(1,.5);
   \draw [xshift=9cm,yshift=2cm] (1,.5)node{}-- (1.5,0)node{} (1,-.5)node{}--(.5,0)node{}
   [xshift=3cm] (1,-.5) node{}-- (1.5,0) node{} (1,.5)node{}--(.5,0)node{};
   \begin{scope}[yshift=-3mm]
     \draw [->] (1.3cm,0cm)--(1.7cm,0cm);
     \draw[->] (1.3cm,2cm)--(1.7cm,2cm);
     \draw[->](8.3cm,0cm)--(8.7cm,0cm);
     \draw[->](8.3cm,2cm)--(8.7cm,2cm);
   \end{scope}
 \end{tikzpicture} 
 \caption{Transforming a member of $\C(H)$
   into  a $2$-factor of  $G$
  (Solid edges represent edges in a member of $\C(H)$
   or a  $2$-factor of  $G$.)}
 \label{fig:cycle}
\end{figure}
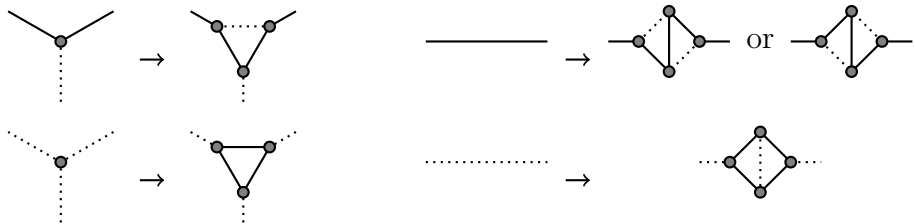

\begin{THM}
  Every claw-free cubic $n$-vertex graph  with no cutedge
  has more than $2^{n/12}$ perfect matchings.
\end{THM}
\begin{proof}
 Let $G$ be a claw-free cubic $n$-vertex graph with no cutedge.
 We may assume that $G$ is connected.
  If $G$ is isomorphic to $K_4$,
  then the claim is clearly true.
  If $G$ is a ring of diamonds,
  then $G$ has  $2^{n/4}+1$ perfect matchings.
 Thus we may assume that $G$ is obtained from a 2-edge-connected
  cubic multigraph $H$ by replacing some edges of $H$ with strings of diamonds
  and replacing each vertex of $H$ with a triangle.
  
  Let $k=|V(H)|$. In other words, $3k$ is the number of vertices not
  in a diamond of $G$.

  Suppose that $k\ge n/6$.
  Since $H$ has $3k/2$ edges,
  the cycle space of $H$ has dimension $3k/2-k+1=k/2+1$
  and therefore $\lvert \C(H)\rvert =2^{k/2+1}$.
  To obtain a $2$-factor from  $C\in \C(H)$,
  we transform $C$ into a member $C'\in \C(G)$
  so that it meets  all $3$ vertices  of $G$ corresponding to $v$
  for each vertex $v$ of $H$ incident with $C$
  as well as it meets all the vertices in each diamond
  that corresponds to an edge in $C$.
 Then  for each vertex $w$ of $G$ unused yet in $C'$,
  we add a cycle of length $3$ or $4$ depending on whether the vertex
  is in a  diamond;
  see Figure~\ref{fig:cycle}.
  Then this is a $2$-factor of $G$ because it meets every vertex of $G$.
  Since the complement of the edge-set of a  2-factor  is a perfect
  matching,
  we conclude that $G$ has at least $2^{k/2+1}\ge 2^{n/12+1}$ perfect matchings.

  Now let us assume that $k<n/6$.
  We know that $G$ has $(n-3k)/4$ diamonds.
  The \emph{length} of an edge $e$ of $H$ is the number of diamonds
  in the string of diamonds replaced with $e$.
  (If the edge $e$ is not replaced with a string of diamonds,
  then the length of $e$ is 0.)
  
  Edmonds' characterization of the perfect matching polytope~\cite{Edmonds1965a}
  implies that
  there exist a positive integer $t$ depending on $H$
  and a list of $3t$ perfect matchings $M_1$, $M_2$,
  $\ldots$, $M_{3t}$ in $H$ such that
  every edge of $H$ is in exactly $t$ of the perfect matchings.
  (In other words, $H$ is fractionally $3$-edge-colorable.)
  By taking complements, we have a list of $3t$ $2$-factors of $H$
 such that each edge of $H$ is in exactly $2t$ of the $2$-factors in the list.
  Since $G$ has $(n-3k)/4$ diamonds,
  the sum of the length of all edges of $H$ is $(n-3k)/4$.
  Therefore there exists a $2$-factor $C$ of $H$
  whose length is at least $\frac{n-3k}{4} \frac 23=(n-3k)/6$.
  
  We claim that $G$ has at least $2^{(n-3k)/6}$ 2-factors
  corresponding to $C$.
  For each diamond in the string replacing an edge $e$ of $C$, 
  there are two ways to route cycles of $C$ through the diamond, see
  Figure~\ref{fig:cycle}. 
  Since $C$ passes through at least $(n-3k)/6$ diamonds,
  $G$ has at least $2^{(n-3k)/6}$ 2-factors.
  Since $k<n/6$,
  $G$ has more than $2^{n/12}$ 2-factors.
  Thus $G$ has more than $2^{n/12}$ perfect matchings.
\end{proof}
 We remark that
every 3-edge-connected claw-free cubic $n$-vertex graph $G$
has exactly $2^{n/6+1}$ perfect matchings, unless $G$ is isomorphic to $K_4$.
That is because $G$ has no
diamonds and so, from the idea of the above proof,  there is a one-to-one correspondence
between the set of all 2-factors of $G$
and the cycle space of a multigraph $H$ obtained by contracting each triangle of $G$.

\end{document}